\documentclass[10pt]{amsart}
\usepackage{amssymb}
\usepackage{bm}
\usepackage{graphicx}
\usepackage[centertags]{amsmath}
\usepackage{amsfonts}
\usepackage{amsthm}
\newtheorem{thm}{Theorem}
\newtheorem{coro}[thm]{Corollary}

\newtheorem{defn}[thm]{Definition}
\newtheorem{pro}{Problem}
\theoremstyle{definition}

\newcommand{\rr}{\mathbb{R}}
\newcommand{\nn}{\mathbb{N}}
\newcommand{\ee}{\varepsilon}

\newcommand{\om}{\omega}
\newcommand{\ct}{2^{<\mathbb{N}}}

\begin{document}

\title{Banach spaces and Ramsey Theory: some open problems}
\author{Pandelis Dodos, Jordi Lopez-Abad and Stevo Todorcevic}

\address{Department of Mathematics, University of Athens, Panepistimiopolis 157 84, Athens, Greece.}
\email{pdodos@math.uoa.gr}

\address{Instituto de Ciencias Matematicas (ICMAT), CSIC-UAM-UC3M-UCM, Madrid, Spain.}
\email{abad@logique.jussieu.fr}

\address{Universit\'{e} Denis Diderot - Paris 7, C.N.R.S., UMR 7056, 2 place Jussieu
- Case 7012, 72521 Paris Cedex 05, France}
\address{and}
\address{Department of Mathematics, University of Toronto, Toronto, Canada, M5S 2E4.}
\email{stevo@math.toronto.edu}


\maketitle


\begin{abstract}
We discuss some open problems in the Geometry of Banach spaces having Ramsey-theoretic
flavor. The problems are exposed together with well known results related to them.
\end{abstract}


\section{Introduction}

One of the most illuminating imports in Banach Space Theory is the field
of Ramsey Theory. There have been some remarkable achievements using combinatorial
techniques; see, for instance, \cite{AGR,AT,G3,Od} and the references therein.
The purpose of this paper is to collect a number of \textit{open problems}
in the Geometry of Banach spaces having Ramsey flavor. We have not tried
to be encyclopedic in our selection. Instead, we have chosen to discuss the
following list of open problems which shows, in particular, that Ramsey Theory
cuts across all parts of modern Banach Space Theory.
\begin{enumerate}
\item[$\bullet$] \textit{Distance from the cube}.
\item[$\bullet$] \textit{Elton unconditionality constant}.
\item[$\bullet$] \textit{Rosenthal basis problem}.
\item[$\bullet$] \textit{Unique spreading model problem}.
\item[$\bullet$] \textit{Block homogeneous basis problem}.
\item[$\bullet$] \textit{Fixing properties of operators on $C[0,1]$}.
\item[$\bullet$] \textit{Bounded distortion}.
\item[$\bullet$] \textit{Converse Aharoni problem}.
\item[$\bullet$] \textit{Separable quotient problem}.
\item[$\bullet$] \textit{Quotients with long Schauder bases and biorthogonal systems}.
\item[$\bullet$] \textit{Rolewicz problem on support sets}.
\item[$\bullet$] \textit{Unconditional basic sequences in non-separable spaces}.
\end{enumerate}
In several cases, a problem naturally leads to a number of related questions
which are interesting on their own. We discuss these issues in detail.

Our Banach space theoretic notation and terminology is standard and follows
\cite{LT}. By $\nn=\{0,1,2,...\}$ we shall denote the natural numbers. If
$S$ is a set, then the cardinality of $S$ will be denoted by $|S|$. Finally,
by $[\nn]^{\infty}$ we denote the set of all infinite subsets of $\nn$.


\section{Distance from the cube}

Recall that the \textit{Banach-Mazur distance} between two isomorphic Banach spaces
$X$ and $Y$ (not necessarily infinite-dimensional) is defined by
\[ d_{BM}(X,Y)=\inf\big\{\|T\|\cdot\|T^{-1}\|: T:X\to Y \text{ is an isomorphism}\big\}.\]
For every $n\in\nn$ define
\[ R^n_\infty=\max\big\{ d_{BM}(X,\ell^n_\infty): \mathrm{dim}(X)=n\big\}.\]
\begin{pro} \label{p1}
\emph{(Distance from the cube -- A. Pe{\l}czy\'{n}ski \cite{Pe2})}
Determine the asymptotic behavior of $R^n_\infty$ as $n$ tends to infinity.
\end{pro}
The best known lower bound for $R^n_\infty$ is due to S. Szarek \cite{Szarek}:
\[ c\sqrt{n}\log n\leqslant R^n_\infty \]
where $c$ is an absolute constant. The problem of finding good upper
estimates for $R^n_\infty$ turned out to be a very challenging problem and
has attracted the attention of several researchers; see \cite{BS,ST}.
The best known upper bound is due to A. Giannopoulos \cite{Gia}:
\[ R^n_\infty \leqslant C n^{5/6}\]
where $C$ is an absolute constant.

This problem appears in our list because all proofs establishing a non-trivial upper bound for $R^n_\infty$
use some variant of the \textit{Sauer-Shelah Lemma}. This basic combinatorial fact asserts that if $S$ is a
non-empty finite set and $\mathcal{F}$ is a ``large" family of subsets of $S$, then there exists a ``large"
subset $A$ of $S$ such that the \textit{trace $\mathcal{F}[A]$ of $\mathcal{F}$ on $A$}, i.e. the family
$\mathcal{F}[A]=\{F\cap A: F\in\mathcal{F}\}$, is the powerset of $A$.
\begin{thm} \label{t1}
\emph{(Sauer-Shelah Lemma \cite{Sauer,Sh1})} Let $S$ be a non-empty finite set, $k\in\nn$ with $1\leqslant
k\leqslant |S|$ and $\mathcal{F}$ be a family of subsets of $S$ such that
\[ |\mathcal{F}| > \sum_{i=0}^{k-1} \binom{|S|}{i}.\]
Then there exists a subset $A$ of $S$ with $|A|=k$ such that the trace $\mathcal{F}[A]$ of $\mathcal{F}$ on
$A$ is the powerset of $A$.
\end{thm}
The Sauer-Shelah Lemma turned out to be a very important tool in the ``local theory" of Banach spaces and has
found significant applications. It was first used in Banach Space Theory by J. Elton \cite{Elt2}.


\section{Elton unconditionality constant}

The Sauer-Shelah Lemma has an infinite version. To state it let us recall that a family $\mathcal{F}$ of finite
subsets of $\nn$ is said to be \textit{compact} (respectively, \textit{hereditary}) if $\mathcal{F}$ is a closed
subset of $2^\nn$ (respectively, for every $F\in\mathcal{F}$ and every $G\subseteq F$ we have that $G\in\mathcal{F}$).
\begin{thm} \label{t2}
\emph{(J. Elton -- see \cite{AGR})}
Let $\mathcal{F}$ be a compact family of finite subsets of $\nn$.
Then there exists $L\in[\nn]^\infty$ such that
$\mathcal{F}[L]$ is a hereditary family.
\end{thm}
This result can be seen as the discrete analogue of a heavily investigated
phenomenon in Banach Space Theory known as \textit{partial unconditionality}.

Recall that a normalized sequence $(x_n)$ is a Banach space $X$ is
said to be \textit{unconditional} if there exists a constant $C\geqslant 1$
such that for every $k\in\nn$, every $a_0,..., a_k\in \rr$ and
every $G\subseteq \{0,...,k\}$ we have
\begin{equation} \label{e1}
\big\| \sum_{n\in G} a_n x_n \big\| \leqslant C \big\| \sum_{n=0}^k a_n x_n \big\|.
\end{equation}
A classical discovery due to B. Maurey and H. P. Rosenthal \cite{MR}
provides an example of a normalized weakly null sequence $(e_n)$ in
a Banach space $E$ with no unconditional subsequence. A surprising
fact is, however, that for \textit{every} normalized weakly
null sequence $(x_n)$ there exist a subsequence $(y_n)$ of $(x_n)$
and a constant $C\geqslant 1$ such that inequality (\ref{e1}) is valid
provided that we restrict our attention to certain sequences of
coefficients (or if we project to certain finite subsets of $\nn$);
i.e. the sequence $(y_n)$ is \textit{partially} unconditional.

By now several partial unconditionality results have been obtained;
see, for instance, \cite{AGR,AT,DOSZ,LTo,Od} and the references therein.
The first (and undoubtedly the most influential) one was discovered by
J. Elton. To state it let us recall the following notion.
\begin{defn} \label{d3}
Let $0<\delta\leqslant 1$ and $K\geqslant 1$. A normalized basic sequence $(x_n)$ is said to be \emph{$\delta$-near unconditional
with constant $K$} if the basis constant of $(x_n)$ is at most $K$ and is such that for every $k\in\nn$, every $a_0,...,a_k\in\rr$
and every $G\subseteq\big\{ i\in\{0,...,k\}: \delta\leqslant |a_i|\leqslant 1\big\}$ we have
\begin{equation} \label{e2}
\big\| \sum_{n\in G} a_n x_n \big\| \leqslant K \big\| \sum_{n=0}^k a_n x_n \big\|.
\end{equation}
\end{defn}
We can now state Elton's Theorem.
\begin{thm} \label{t4}
\emph{(J. Elton \cite{Elt1})}
Let $0<\delta\leqslant 1$. Then there exists a constant $K\geqslant 1$, possibly depending on $\delta$, such that
for every normalized weakly null sequence $(e_n)$ in a Banach space $E$ there exists a subsequence $(x_n)$ of $(e_n)$
such that $(x_n)$ is $\delta$-near unconditional with constant $K$.
\end{thm}
For each $0<\delta\leqslant 1$ let $K(\delta)$ be the infimum of the set
of real numbers $K$ such that every normalized weakly null sequence has
a $\delta$-near unconditional subsequence with constant $K$.
\begin{pro} \label{p2}
\emph{(Elton unconditionality constant -- S. J. Dilworth, E. Odell, Th. Schlumprecht and A. Zs\'{a}k
\cite{DOSZ})} Does there exist $M\geqslant 1$ such that $K(\delta)\leqslant M$ for every $0<\delta\leqslant
1$?
\end{pro}
It is known (see \cite{DKK,Od}) that there exists an absolute constant
$C$ such that
\[ K(\delta)\leqslant C\log(1/\delta)\]
for every $0<\delta\leqslant 1$. It is also known (see \cite{DOSZ}) that
$5/4\leqslant \sup\{K(\delta):0<\delta\leqslant 1\}$.


\section{Rosenthal basis problem}

A normalized basic sequence $(x_n)$ in a Banach space $X$ is said to
be \textit{perfectly homogeneous} if every normalized block sequence
$(v_n)$ of $(x_n)$ is equivalent to $(x_n)$. It is easy to see that
the standard unit vector bases of the classical sequence spaces $c_0$
and $\ell_p$ $(1\leqslant p<+\infty)$ are perfectly homogeneous. Remarkably,
the converse is also true.
\begin{thm} \label{t5}
\emph{(M. Zippin \cite{Zi})}
Every perfectly homogeneous basic sequence $(x_n)$ is equivalent to the
standard unit vector basis of $c_0$ or $\ell_p$ for some $1\leqslant p<+\infty$.
\end{thm}
There is a slightly weaker notion: a normalized basic sequence $(x_n)$ in a Banach
space $X$ is said to be a \textit{Rosenthal basis} if every normalized block
sequence $(v_n)$ of $(x_n)$ has a subsequence which is equivalent to $(x_n)$.
\begin{pro} \label{p3}
\emph{(Rosenthal basis problem -- H. P. Rosenthal \cite{FPR})}
Is it true that every Rosenthal basis $(x_n)$ is equivalent to the standard
unit vector basis of $c_0$ or $\ell_p$ for some $1\leqslant p<+\infty$?
\end{pro}
It is known (see \cite{FPR}) that the problem has an affirmative answer provided
that there exists a constant $C\geqslant 1$ such that every normalized block
sequence $(v_n)$ of $(x_n)$ has a subsequence which is $C$-equivalent to $(x_n)$.


\section{Unique spreading model problem}

Let $(x_n)$ be a normalized basic sequence in a Banach space $X$ and $(y_n)$ be a sequence
in a Banach space $Y$. The sequence $(x_n)$ is said to \textit{generate the sequence $(y_n)$
as spreading model} if there exists a sequence $(\ee_n)$ of positive reals with $\ee_n\downarrow 0$
such that for every $n\in\nn$, every $a_0,...,a_n$ in $[-1,1]$ and every
$n\leqslant k_0< ... < k_n$ in $\nn$ we have
\[ (1+\ee_n)^{-1} \big\| \sum_{i=0}^n a_i y_i\big\| \leqslant
\big\| \sum_{i=0}^n a_i x_{k_i} \big\| \leqslant
(1+\ee_n) \big\| \sum_{i=0}^n a_i y_i\big\|.\]
It is well-known (see \cite{BrS1,BrS2}) that every normalized basic
sequence has a subsequence generating a spreading model. Actually, this
fact is one of the earliest applications of Ramsey Theory in the
Geometry of Banach spaces.

For every Banach space $X$ let $\mathrm{SP_w}(X)$ be the class of all
spreading models generated by a normalized weakly null sequence in $X$
and denote by $\sim$ the usual equivalence relation of equivalence of
basic sequences.
\begin{pro} \label{p4}
\emph{(Unique spreading model problem -- S. A. Argyros \cite{AOST})}
Let $X$ be a reflexive Banach space such that $|\mathrm{SP_w}(X)/\sim|=1$ (i.e. all
normalized weakly null sequences in $X$ generate a unique spreading model
up to equivalence). Is it true that the unique spreading model is equivalent
to the standard unit vector basis of $c_0$ or $\ell_p$ for some
$1\leqslant p<+\infty$?
\end{pro}
It is known (see \cite{DOS}) that the reflexivity assumption is essential.
It is also known (and easy to see) that an affirmative answer to
Problem \ref{p4} yields an affirmative answer to Problem \ref{p3}.

The following stronger version of Problem \ref{p4} has been also considered.
\begin{pro} \label{p5}
\emph{(Countable spreading models problem -- S. J. Dilworth, E. Odell and B. Sari \cite{DOS})}
Let $X$ be a reflexive Banach space such that $|\mathrm{SP_w}(X)/\sim|\leqslant \aleph_0$.
Does there exist a normalized weakly null sequence in $X$ generating a spreading model
equivalent to the standard unit vector basis of $c_0$ or  $\ell_p$ for some
$1\leqslant p<+\infty$?
\end{pro}
Some partial answers have been obtained. For instance, it is known (see \cite{DOS})
that if $X$ is reflexive, $|\mathrm{SP_w}(X)/\sim|\leqslant \aleph_0$ and
$|\mathrm{SP_w}(Y^*)/\sim|\leqslant \aleph_0$ for every infinite-dimensional subspace
$Y$ of $X$, then there exists a normalized weakly null sequence in $X$ generating a spreading model
equivalent to the standard unit vector basis of $c_0$ or $\ell_p$ for some
$1\leqslant p<+\infty$.


\section{Block homogeneous basis problem}

A normalized Schauder basis $(x_n)$ of a Banach space $X$ is said to be
\textit{block homogeneous} if every normalized block sequence $(v_n)$
of $(x_n)$ spans a subspace isomorphic to $X$.
\begin{pro} \label{p6}
\emph{(Block homogeneous basis problem -- see \cite{FR})}
Let $X$ be a Banach space with a normalized block homogeneous Schauder
basis. Is it true that $X$ is isomorphic to $c_0$ or $\ell_p$ for some
$1\leqslant p<+\infty$?
\end{pro}
Notice that an affirmative answer to the above problem would considerably
strengthen Theorem \ref{t5} and, combined with the deep results of
A. Szankowski \cite{Szankowski}, would yield an alternative proof
to the \textit{Homogeneous Banach Space Problem} that avoids the use of
W. T. Gowers' dichotomy \cite{G2}.

Although it is known that every separable non-Hilbertian Banach space contains
two non-isomorphic subspaces, we point out that there are several
related open problems. The following one is probably the most natural.
\begin{pro} \label{p7}
\emph{(G. Godefroy \cite{Go})}
Let $X$ be a separable non-Hilbertian Banach space. Is it true that there
exists a sequence $(X_n)$ of subspaces of $X$ which are pairwise non-isomorphic?
\end{pro}
It is likely that Problem \ref{p7} has an affirmative answer. Indeed, several
partial results have been obtained (see \cite{FR} and the references therein).

We would like to isolate an instance of Problem \ref{p7} exposing, in particular,
its connection with the \textit{Basis Problem}, a famous problem in the Geometry
of Banach spaces that has received considerable attention during the 1970s and 1980s.
Recall that, by the results in \cite{MP,Szankowski}, if a separable Banach space
$X$ fails to have type $2-\ee$ (or fails to have co-type $2+\ee$) for some
$\ee>0$, then $X$ contains a subspace without a Schauder basis (actually,
$X$ contains a subspace without the compact approximation property).
\begin{pro} \label{p8}
Let $X$ be a separable Banach space which fails to have type $2-\ee$ (or fails
to have co-type $2+\ee$) for some $\ee>0$. Is it true that $X$ contains
continuum many pairwise non-isomorphic subspaces?
\end{pro}
It is known (see \cite{Ani}) that the above problem has an affirmative answer
if $X$ is a separable weak Hilbert space not isomorphic to $\ell_2$.


\section{Fixing properties of operators on $C[0,1]$}

A central open problem in Banach Space Theory is the classification, up to isomorphism,
of all complemented subspaces of the classical function space $C[0,1]$.
\begin{pro} \label{csc}
\emph{(Complemented Subspace Problem for $C[0,1]$ -- see \cite{BP-new,LW-new,Pe1})}
Is it true that every complemented subspace of $C[0,1]$ is isomorphic to a $C(K)$ space
for some closed subset $K$ of $[0,1]$.
\end{pro}
The \textit{Complemented Subspace Problem} is discussed in detail in
\cite{Ro2}. The strongest result to date is due to H. P. Rosenthal
\cite{Ro1} and asserts that every complemented subspace of $C[0,1]$ with
non-separable dual is isomorphic to $C[0,1]$. Actually, this fact is a
consequence of the following surprising discovery concerning general
operators on $C[0,1]$.
\begin{thm} \label{t6}
\emph{(H. P. Rosenthal \cite{Ro1})}
Let $T:C[0,1]\to C[0,1]$ be a bounded linear operator such that its dual $T^*$
has non-separable range. Then $T$ fixes a copy of $C[0,1]$; that is, there
exists a subspace $Y$ of $C[0,1]$ which is isomorphic to $C[0,1]$ and is such
that $T|_Y$ is an isomorphic embedding.
\end{thm}
The following problem, asked in the 1970s, is motivated by the
\textit{Complemented Subspace Problem} and Theorem \ref{t6}.
\begin{pro} \label{p9}
\emph{(Fixing properties of operators on $C[0,1]$ -- H. P. Rosenthal \cite{Ro2})}
Let $T:C[0,1]\to C[0,1]$ be a bounded linear operator and $X$ be a separable Banach
space not containing a copy $c_0$. Assume that the operator $T$ fixes a copy of $X$.
Is it true that $T$ fixes a copy of $C[0,1]$? Equivalently, is it true that the dual
operator $T^*$ of $T$ has non-separable range?
\end{pro}
Notice that an affirmative answer to Problem \ref{p9} yields that every
complemented subspace $X$ of $C[0,1]$ with separable dual is hereditarily $c_0$
(that is, every infinite-dimensional subspace of $X$ contains a copy of $c_0$).
This property is predicted by the \textit{Complemented Subspace Problem}
(see \cite{PS}) though it has not been decided yet. It is known, however,
that every complemented subspace of $C[0,1]$ contains a copy $c_0$ (see \cite{Pe1}).

There are some partial positive answers to Problem \ref{p9}.
\begin{thm} \label{tB-new}
\emph{(J. Bourgain \cite{Bou})}
If an operator $T:C[0,1]\to C[0,1]$ fixes a copy of a co-type space, then $T$ fixes a copy of $C[0,1]$
\end{thm}
\begin{thm} \label{tG-new}
\emph{(I. Gasparis \cite{Ga})}
If an operator $T:C[0,1]\to C[0,1]$ fixes a copy of an asymptotic $\ell_1$ space, then $T$ fixes a copy of $C[0,1]$.
\end{thm}
Also, recently, a trichotomy was obtained characterizing, by means of fixing properties, the
class of operators on general separable Banach spaces whose dual has non-separable range.
To state it let us say that a sequence $(x_t)_{t\in\ct}$ in a Banach space $X$ is said
to be a \textit{tree basis} provided that the following are satisfied.
\begin{enumerate}
\item[(1)] If $(t_n)$ is the canonical enumeration\footnote{This means that for every $n,m\in\nn$
we have $n<m$ if and only if either $|t_n|<|t_m|$ or $|t_n|=|t_m|$ and $t_n<_{\mathrm{lex}} t_m$.}
of $\ct$, then $(x_{t_n})$ is a seminormalized basic sequence.
\item[(2)] For every infinite antichain $A$ of $\ct$ the sequence $(x_t)_{t\in A}$ is weakly null.
\item[(3)] For every $\sigma\in 2^\nn$ the sequence $(x_{\sigma|n})$ is weak* convergent to an element
$x^{**}_\sigma\in X^{**}\setminus X$. Moreover, if $\sigma, \tau\in 2^\nn$ with $\sigma\neq \tau$,
then $x^{**}_\sigma\neq x^{**}_\tau$.
\end{enumerate}
The archetypical example of such a sequence is the standard unit vector basis
of James tree space (see \cite{Ja-new}).
\begin{thm} \label{tD-new}
\emph{(P. Dodos \cite{D})}
Let $X$ and $Y$ be separable Banach spaces. Also let $T:X\to Y$ be an operator.
\begin{enumerate}
\item[(a)] If $X$ does not contain a copy of $\ell_1$, then $T^*$ has non-separable range
if and only if there exists a sequence $(x_t)_{t\in\ct}$ in $X$ such that both $(x_t)_{t\in\ct}$
and $(T(x_t))_{t\in\ct}$ are tree bases.
\item[(b)] If $X$ contains a copy of $\ell_1$, then $T^*$ has non-separable range if and only if
either the operator $T$ fixes a copy of $\ell_1$, or there exists a bounded sequence
$(x_t)_{t\in\ct}$ in $X$ such that its image $(T(x_t))_{t\in\ct}$ is a tree basis.
\end{enumerate}
\end{thm}
The proofs of the aforementioned results have a strong combinatorial flavor. In particular,
in \cite{Ga} heavy use is made of the infinite-dimensional Ramsey Theorem for
Borel partitions (see \cite{GP}) while the main result in \cite{D} is based on
the deep Halpern-L\"{a}uchli Theorem and its consequences (see \cite{HL,Ka,To2}).


\section{Bounded distortion}

Let $(X, \|\cdot\|)$ be an infinite-dimensional Banach space and $\lambda>1$.
The space $X$ is said to be \textit{$\lambda$-distortable} if there exists an
equivalent norm $|\cdot |$ on $X$ such that for every infinite-dimensional
subspace $Y$ of $X$ we have
\[ \sup\Big\{ \frac{|y|}{|z|}: y, z \in Y \text{ and } \|y\|=\|z\|=1\Big\} \geqslant \lambda.\]
The space $X$ is said to be \textit{distortable} if it is $\lambda$-distortable
for some $\lambda>1$. A distortable Banach space $X$ is said to be
\textit{arbitrarily distortable} if it is $\lambda$-distortable for
every $\lambda>1$; otherwise, $X$ is said to be \textit{boundedly distortable}.

The first result on distortion is due to R. C. James \cite{Ja}: \textit{$c_0$
and $\ell_1$ are not distortable}. V. D. Milman \cite{Mi} proved that
a non-distortable Banach space must contain almost isometric copies of
$c_0$ or $\ell_p$ for some $1\leqslant p <+\infty$. The first example of a
distortable Banach space was Tsirelson's space $T$ \cite{Ts}. The natural
problem whether the spaces $\ell_p$ $(1<p<+\infty)$ are distortable
became known as the \textit{Distortion Problem} and it was solved by
the following famous result.
\begin{thm} \label{t7}
\emph{(E. Odell and Th. Schlumprecht \cite{OS})}
For any $1<p<+\infty$ the space $\ell_p$ is arbitrarily distortable.
\end{thm}
The \textit{Distortion Problem} is discussed in detail in \cite{OS-new}.
A surprising phenomenon observed so far is that ``most" distortable
Banach spaces are arbitrarily distortable.
\begin{pro} \label{p10}
\emph{(Bounded distortion -- see \cite{OS-new})} Does there exist a boundedly distortable
Banach space?
\end{pro}
Tsirelson's space is a candidate for being boundedly distortable.
Indeed, we have already mentioned that $T$ is distortable;
actually, $T$ is $(2-\ee)$-distortable for every $\ee>0$ (see \cite{BL}).
But it not known whether $T$ is arbitrarily distortable.
\begin{pro} \label{p11}
\emph{(Distortion of Tsirelson's space -- H. P. Rosenthal)}
Is $T$ (or its dual) boundedly distortable?
\end{pro}
The \textit{Distortion Problem} has a non-linear analogue which
is also closely related to Ramsey Theory. Recall that a real-valued
Lipschitz map $f$ defined on the unit sphere $S_X$ of an infinite-dimensional
Banach space $X$ is said to be \textit{oscillation stable} if for every
$\ee>0$ there exists an infinite-dimensional subspace $Y$ of $X$ such that
the oscillation of $f$ on $S_Y$ is at most $\ee$, i.e.
\[ \sup\{ |f(y)-f(z)|: y, z\in Y \text{ and } \|y\|=\|z\|=1\}\leqslant \ee.\]
Notice that if a Banach space $X$ has the property that every real-valued
Lipschitz map $f$ on $S_X$ is oscillation stable, then $X$ is not distortable.
This is, however, a stronger property and is possessed only by hereditarily
$c_0$ spaces (see \cite{BL}). That the space $c_0$ has this property follows
by the following remarkable result.
\begin{thm} \label{t8}
\emph{(W. T. Gowers \cite{G1})} Every real-valued uniformly continuous
map on the unit sphere of $c_0$ is oscillation stable.
\end{thm}
All known proofs of Theorem \ref{t8} use methods of topological dynamics
on $\beta\nn$. It would be very interesting to have a purely combinatorial
proof avoiding the use of ultrafilters.


\section{Converse Aharoni problem}

The space $c_0$ is admittedly ``small" in the linear category: \textit{if
a Banach space $X$ is isomorphic to a subspace of $c_0$, then $X$ contains
a copy of $c_0$}. The situation, however, is dramatically different in
the non-linear category.
\begin{thm} \label{t9}
\emph{(I. Aharoni \cite{Ah})} Every separable metric space is Lipschitz
isomorphic to a subset of $c_0$.
\end{thm}
The best constant of the Lipschitz embedding has been computed (see \cite{KL}).
Theorem \ref{t9} can be rephrased to say that $c_0$ is universal in the
non-linear category. But it is not known whether $c_0$ is the
\textit{minimal} Banach space with this property.
\begin{pro} \label{p12}
\emph{(Converse Aharoni problem -- see \cite{BL,Kalton})}
Let $X$ be a separable Banach space with the property that every separable
Banach space is Lipschitz isomorphic to a subset of $X$. Is it
true that $X$ contains a copy of $c_0$?
\end{pro}
It is known (see \cite{GKL}) that the class of subspaces of $c_0$ is stable
under Lipschitz isomorphisms (moreover, if a Banach space $X$ is Lipschitz
isomorphic to $c_0$, then $X$ is linearly isomorphic to $c_0$). It is
possible that Theorem \ref{t8} could be used to attack Problem \ref{p12}.


\section{Separable quotient problem}

In the last four sections we will present some problems about general Banach spaces
where we do not assume separability. It turns out that this area has its own Ramsey
Theory to which it is related, the rich area of Ramsey Theory of uncountable structure
with a strong set-theoretic flavor (see, for example, \cite{To0,To5,To6}).

Recall the very old result due to S. Banach and S. Mazur saying that every
infinite-dimensional Banach space has an infinite basic sequence. It turns out
that the dual version of this \textit{Basic Sequence Problem} is still widely open.
\begin{pro} \label{p13}
\emph{(S. Banach \cite{B} -- A. Pe{\l}czy\'{n}ski \cite{Pe0})}
Does every infinite-dimensional Banach space $X$ has an infinite-dimensional
quotient with a Schauder basis?
\end{pro}
The following result represents the first major progress made on this problem.
\begin{thm} \label{t10}
\emph{(W. B. Johnson and H. P. Rosenthal \cite{JR})}
Every infinite-dimensional separable space $X$, and in fact every infinite-dimensional
space of density $<\mathfrak{b},$\footnote{Recall that $\mathfrak{b}$ is the minimal
cardinality of a subset of $\nn^\nn$ which is unbounded in the ordering of eventual dominance.}
has an infinite-dimensional quotient with a Schauder basis.
\end{thm}
So the problem of S. Banach and A. Pe{\l}czy\'{n}ski now becomes the following.
\begin{pro} \label{p14}
\emph{(Separable quotient problem)} Does every infinite-dimensional Banach
space have a non-trivial infinite-dimensional separable quotient?
\end{pro}
There are several partial results worth mentioning (see \cite{Mu} for a rather
complete historical review of the early work on this problem).
\begin{thm} \label{t11}
\emph{(W. B. Johnson and H. P. Rosenthal \cite{JR} -- J. N. Hagler and W. B. Johnson \cite{HJ})}
If $X^*$ contains an infinite-dimensional subspace with separable dual,
then $X$ has a non-trivial separable quotient. If $X^*$ has an unconditional
basic sequence, then $X$ has a quotient with an unconditional basis.
\end{thm}
\begin{thm} \label{t12}
\emph{(S. A. Argyros, P. Dodos and V. Kanellopoulos \cite{ADK1,ADK2})}
Every dual Banach space and every representable Banach space has a separable
non-trivial quotient.
\end{thm}
Concerning the two sufficient conditions given in Theorem \ref{t11} above,
we mention that W. T. Gowers \cite{G-new} has constructed a \emph{separable}
dual space which contains neither a copy of $\ell_1$ nor an infinite-dimensional
subspace with a separable dual.


\section{Quotients with long Schauder bases and biorthogonal systems}

Note that Problem \ref{p13} has also the following natural variant.
\begin{pro} \label{p15}
Given an infinite-dimensional Banach space $X$, what is the longest Schauder basis
a quotient of $X$ can have?
\end{pro}
Note that this in particular asks if every non-separable Banach space has an
uncountable biorthogonal system, a problem that has been first explicitly
addressed by W. J. Davis and W. B. Johnson \cite{DJ}. Today we know that the
problem is really about additional set-theoretic methods that can be relevant
in building uncountable biorthogonal systems and quotient maps. The first
and today well-known examples showing this were constructed by K. Kunen (see
\cite{Ku1,Ku2,N}) assuming CH and S. Shelah \cite{Sh2} assuming $\diamondsuit$.
Recently a systematic study of generic Banach spaces is given in \cite{LTo2}
and many of these examples are related to the uncountable biorthogonal system
problem. The following old example (see \cite[Chapter 1]{To0}) is not that
well-known and we mention it here because of its relevance to the discussion below.
\begin{thm} \label{t13}
\emph{(S. Todorcevic \cite{To0})}
If $\mathfrak{b}=\aleph_1$, then there is an Asplund $C(K)$ space without uncountable
biorthogonal systems.
\end{thm}
The second version of Problem \ref{p13} (Problem \ref{p15}) was analyzed in \cite{To1}
where the problem was connected with the following well-known set-theoretic and
Ramsey-theoretic principle.
\medskip

\noindent \textbf{P-ideal dichotomy (PID)}: For every P-ideal\footnote{Recall that
an ideal $\mathcal{I}$ is a P-ideal if for every sequence $(X_n)$ of elements of
$\mathcal{I}$ there is $Y$ such that $X_n\setminus Y$ is finite for all $n$.}
$\mathcal{I}$ consisting of countable subsets of some index-set $S$, either
\begin{enumerate}
\item[(1)] there is an uncountable set $T\subseteq S$ such that every countable
subset of $T$ belongs to $\mathcal{I}$, or
\item[(2)] the set $S$ can be decomposed into countably many subsets that have
no infinite subsets belonging to $\mathcal{I}$.
\end{enumerate}
\begin{thm} \label{t14}
\emph{(S. Todorcevic \cite{To1})} Assume {\rm PID}. Every nonseparable Banach space
$X$ of density $<\mathfrak{m}$ has a quotient with a Schauder basis of length $\omega_1$.
\end{thm}
Recall that $\mathfrak{m}$ is the minimal cardinality of a family $\mathcal{F}$ of
nowhere dense subsets of a compact ccc space $K$ such that $\bigcup\mathcal{F}=K$
(see, for example, \cite{Fr1}).
\begin{coro} \label{c15}
\emph{(S. Todorcevic \cite{To1})}
Assume {\rm PID}. If $\mathfrak{m}>\aleph_1$\footnote{Under {\rm PID}, the statement
$\mathfrak{m}>\aleph_1$ and the statement $\mathfrak{m}=\aleph_2$ are equivalent.
This holds also for other cardinals such as $\mathfrak{b}$ and $\mathfrak{p}$
(see \cite{To3,To4}).}, then every non-separable Banach space has an uncountable
biorthogonal system.
\end{coro}
The reader is referred to \cite{To3} and \cite{To4} which list known applications
of this principle many of which  have a strong Ramsey-theoretic flavor. It turns
out that {\rm PID} is an additional set-theoretic principle which, while independent
from CH, has a power of reducing problems about general sets to problems about reals.
Finding out what is the actual problem about reals that corresponds to the quotient
basis problem or the biorthogonal sequence problem is now in order. The following
is a particular instance of what one may expect from such analysis.
\begin{pro} \label{p16}
Are the following equivalent assuming the P-ideal dichotomy?
\begin{enumerate}
\item[(1)] Every non-separable Banach space has an uncountable biorthogonal system.
\item[(2)] $\mathfrak{b}=\aleph_2$.
\end{enumerate}
\end{pro}


\section{Rolewicz's problem on support sets}

Recall that a \emph{support set} in a Banach space $X$ is a nonempty convex set
$C$ with the property that every point $x$ of $C$ is its support point, i.e. there
is a functional $f$ of $X$ such that
\[ f(x)=\inf\{f(y): y\in C\}<\sup\{f(y): y\in C\}. \]
The existence of a support set in every non-separable Banach space is an old problem
of S. Rolewicz \cite{Rol} who showed that separable Banach spaces do not admit
support sets and who noticed that many of the non-separable Banach spaces
do admit support sets.
\begin{pro} \label{p17}
\emph{(S. Rolewicz \cite{Rol})}
Which Banach spaces admit support sets? Do all non-separable Banach spaces
have support sets?
\end{pro}
It is easily seen that if $X$ has an uncountable biorthogonal system, then $X$ has a support set, so the second part of this
question has a positive answer assuming $\mathrm{PID}$ and $\mathfrak{m}>\aleph_1$. It is also known that some assumptions
are necessary as it is consistence to have a $C(K)$ space of density $\aleph_1$ without support sets (see \cite{LTo2} and
\cite{Kos}). It turns out, however, that finding uncountable biorthogonal systems is a different problem from finding
support sets in Banach spaces.
\begin{thm} \label{t16}
\emph{(S. Todorcevic \cite{To1})} Every $C(K)$ space of density $>\aleph_1$ has a support set.
\end{thm}
Combining this with the following result we see the difference between
the Rolewicz problem and the uncountable biorthogonal sequence problem.
\begin{thm} \label{t17}
\emph{(C. Brech and P. Koszmider \cite{BK})}
It is consistent to have $C(K)$ spaces of density $>\aleph_1$ without
uncountable biorthogonal systems.
\end{thm}
This leaves however open the following interesting variant of Rolewicz's problem.
\begin{pro} \label{p18}
\emph{(S. Todorcevic \cite{To1})} Does every Banach space of density $>\aleph_1$
have a support set?
\end{pro}


\section{Unconditional basic sequences in non-separable spaces}

We now turn our attention to another well-known problem in this area.
\begin{pro} \label{p19}
\emph{(Unconditional basic sequences in non-separable spaces)}
Does every Banach space $X$ of density $>\mathfrak{c}$ have an infinite unconditional basic sequence? Does additional properties,
like reflexivity, make a difference?
\end{pro}
We mention some results related to this problem.
\begin{thm} \label{t18}
\emph{(S. A. Argyros, A. D. Arvanitakis and A. G. Tolias \cite{AAT})}
There is a Banach space $X$ of density $\mathfrak c$ with no infinite unconditional sequences (in fact, $X$ is hereditarily indecomposable)
nor an infinite-dimensional reflexive subspace.
\end{thm}
\begin{thm} \label{t19}
\emph{(S. A. Argyros, J. Lopez-Abad and S. Todorcevic \cite{ALT})}
There exists a reflexive Banach space of density $\aleph_1$ without infinite unconditional sequences.
\end{thm}
So, the unconditional basic sequence problem might have a different answers for different classes of Banach spaces. This has been hinted also
in a recent article by the authors where the following result is proven.
\begin{thm} \label{t20}
\emph{(P. Dodos, J. Lopez-Abad and S. Todorcevic \cite{DLT})}
It is consistent to assume that every Banach space of density $\geqslant \aleph_\om$ has an infinite unconditional sequence.
\end{thm}
\begin{coro} \label{c21}
\emph{(P. Dodos, J. Lopez-Abad and S. Todorcevic \cite{DLT})}
It is consistent to assume that every Banach space of density $\geqslant \aleph_\om$ has a quotient with an unconditional basis.
\end{coro}
There are related questions about the existence of infinite sub-symmetric sequences which seem to call for a different approach.
First of all let us recall the following well-known result of J. Ketonen in a slightly different form.
\begin{thm} \label{t22}
\emph{(J. Ketonen \cite{Ke})} Every Banach space of density at least the first Erd\H{o}s cardinal\footnote[5]{Recall that $\theta$
is said to be an \textit{Erd\H{o}s cardinal} if for every function $f:[\theta]^{<\omega}\to\omega$ there is an infinite set
$\Gamma\subseteq\theta$ such that $f$ is constant on $[\Gamma]^n$ for every $n<\omega$. This is a large cardinal which is larger
than the first inaccessible cardinal and, in particular, much larger than the continuum.} has an infinite sub-symmetric sequence.
\end{thm}
On the other hand no consistent bounds, like the one of \cite{DLT}
mentioned above, is known in this context. So the following problem is wide open.
\begin{pro} \label{p20}
Does every (respectively, every reflexive) Banach space of density $>\mathfrak{c}$
(respectively, of density $>\aleph_1$) has an infinite sub-symmetric basic
sequence?
\end{pro}
While analyzing this problem one is naturally led to the following variation.
\begin{pro} \label{p21}
Is there a non-separable Tsirelson-like space?
\end{pro}
This problem seems related to the existence of a nonseparable analogue of
the Schreier family on which the Tsirelson construction could be based.
This is of course related to the following well-known problem.
\begin{pro} \label{p23}
\emph{(D. H. Fremlin \cite{Fr2})}
Is there a compact hereditary family of finite subsets of $\omega_1$ containing
the singletons and having the property that every finite subset of $\omega_1$
has a subset half of its size belonging to the family $\mathcal{F}$?
\end{pro}

\noindent \textbf{Acknowledgment.} We would like to thank Gilles Godefroy
for his valuable comments and remarks.


\end{document}